*Research Article*

# A Fundamental Criteria to Establish General Formulas of Integrals

**Rania Saadeh**,[1] **Mohammad Abu-Ghuwaleh**,[1] **Ahmad Qazza**,[1] **and Emad Kuffi**[2]

[1]*Department of Mathematics, Faculty of Science, Zarqa University, Zarqa 13110, Jordan*
[2]*College of Engineering, Al-Qadisiyah University, Al Diwaniyah, Iraq*

Correspondence should be addressed to Rania Saadeh; rsaadeh@zu.edu.jo

In this study, new master theorems and general formulas of integrals are presented and implemented to solve some complicated applications in different fields of science. The proposed theorems are considered to be generators of new problems, including difficult integrals with their exact solutions. The results of these problems can be obtained directly without the need for difficult calculations. New criteria for treating improper integrals are presented and illustrated in four interesting examples and some tables to simplify the procedure of using the proposed theorems. The outcomes of this study are compared with those presented by Gradshteyn and Ryzhik in the classical table of integrations. The results in this study are simple and applicable in solving integrals, and some of the well-known theorems in calculating improper integrals are considered simple cases of our research.

## 1. Introduction

During recent decades, many studies on the theory of improper integrals have been conducted in different fields of science, such as physics and engineering [1–8]. Hence, these integrals were very attractive for mathematicians to discover new theorems and techniques for solving them. Many applications need improper integrals to handle, either in the calculations or in expressing the models, especially in engineering, applied mathematical physics, electronics engineering, etc. [9–15]. Some of these integrations can be solved simply, but others need difficult and long computations. A large number of these integrals cannot be solved manually, but they need computer software to be solved. Additionally, numerical methods can be used to solve some improper integrals that cannot be solved by previous methods [16–22].

The process of evaluating improper integrals is not usually based on certain rules or techniques that can be applied directly. Many methods and techniques were established and introduced by mathematicians and physicists to present a closed form for indefinite integrals, the technique of double integrals, series methods, residue theorem, calculus under the integral sign, and other methods that are used to solve improper complex integrals exactly or approximately, see [23–28].

The residue theorem was first established by Cauchy in 1826, which is considered a powerful theorem in complex analysis. However, the applications that can be calculated using the residue theorem to compute integrals on real numbers need many precise constraints that should be satisfied to solve the integrals, including finding appropriate closed contours and determining the poles. Another challenge in the process of applying the residue theorem is the difficulty and efforts in finding solutions for some integrations.

In his published memoirs, Cauchy reached powerful formulas in mathematics using the residue theorem [4]. Researchers consider these formulas essential in treating and solving improper integrals. However, these results are considered simple cases compared to the results that we present in this article. In addition, we mention that the proposed theorems and results in this research are not based on the residue theorem.

One significant accomplishment in the sphere of definite and indefinite integrals is the master theorem of Ramanujan,

which presents new expressions concerning the Milline transform of any continuous function in terms of the analytic Taylor's series [29–34]. It was implemented by Ramanujan and other researchers as a powerful tool in calculating definite and indefinite integrals and in computing infinite series. The obtained results are applicable and effective as Ramanujan's master theorem in handling and generating new formulas of integrals with direct solutions.

In this study, we introduce new theorems to simplify the procedure of computing improper integrals by presenting new theorems with proofs. Each theorem can generate many improper integral formulas that cannot be solved by usual techniques or need much effort and time to be solved. The theorems introduce the solutions directly in a simple finite sum that depends on the obtained integral. The motivation of this work is to generate as many improper integrals and their values as possible to be used in different applications and problems. The obtained results can be implemented to construct new tables of integrations so that researchers can use them in calculations and to check the accuracy of their answers during the discovery of new methods.

The main purpose of this work is to introduce simple new techniques to help researchers, mathematicians, engineers, physicists, etc., to solve some difficult improper integrals that cannot be treated or solved easily, which requires several theorems and much effort to solve by presenting new approaches. This goal is achieved by introducing some master theorems that can be implemented to solve difficult applications. The outcomes can be generalized and introduced in tables to be used and obtain the results of some improper integrals directly.

Within these results, we introduce a closed expression of integrals that can be established by defining a suitable function on which the target application depends. We consider these theorems as a solid tool for unravelling new families of improper integrals and creating many complicated and interesting integrals that can be solved directly based on our new results.

We organize this article as follows. In Section 2, we introduce some illustrative preliminaries and facts concerning analytic and special functions. Master theorems and results are presented in Section 3. Mathematical remarks and several applications are presented in Section 4. Finally, the conclusion of our research is presented in Section 5.

## 2. Preliminaries

In this section, some basic definitions and theorems related to our work are presented and illustrated for later use.

*2.1. Basic Definitions and Theorems*

*Definition 1 see* [8]. Suppose that a function $f$ is analytic in a domain $\Omega \subseteq \mathbb{C}$, where $\mathbb{C}$ is the complex plane. Consider a disc $D \subseteq \Omega$ centered at $z_0$; then, the function $f$ can be expressed in the following series expansion:

$$f(z) = \sum_{n=0}^{\infty} a_n (z - z_0)^n. \tag{1}$$

*Definition 2 see* [9]. Assume that $f$ is an analytic function; then, the Taylor series expansion at any point $x_0$ of $f$ in its domain is given by

$$T(x) = \sum_{n=0}^{\infty} \frac{f^{(n)}(x_0)}{n!} (x - x_0)^n, \tag{2}$$

that converges to $f$ in a neighborhood of $x_0$ point wisely.

*Definition 3 see* [10]. The Cauchy principal value of an infinite integral of a function $f$ is defined by

$$PV \int_0^{\infty} f(x)dx \equiv \lim_{R \longrightarrow \infty} \int_0^R f(x)dx. \tag{3}$$

*2.2. Basic Formulas of Series and Improper Integrations.* In this section, we introduce some series and improper integrals that are needed in our work.

**Lemma 1.** *Let $n \in \mathbb{N}$, then*

$$\frac{1}{\left(1^2 - x^2\right)\left(3^2 - x^2\right) \cdots \left((2n+1)^2 - x^2\right)} = \frac{(-1)^{n+1}}{4^n (2n+1)!} \sum_{k=0}^{n} (-1)^k \binom{2n+1}{k} \frac{2n+1-2k}{x^2 - (2n+1-2k)^2}. \tag{4}$$

*Proof.* To prove equation (4), we define an integral whose solution can be expressed by two different forms: the left side of equation (4) and the right side of the equation.

Let

$$I = \int_0^{\infty} e^{-xt} (\sinh t)^{2n+1} dt, \tag{5}$$

where $p > 0, x > 2n + 1, n \in \mathbb{N}$.

Taking the indefinite integral:

$$J = x^2 \int e^{-xt} (\sinh t)^{2n+1} dt. \tag{6}$$

Applying integration by parts on equation (6) twice, we obtain a reduction formula as follows:

$$\begin{aligned} J = &-e^{-xt} (\sinh t)^{2n+1} - (2n+1) e^{-xt} \sinh^{2n} t \cosh t \\ &+ (2n+1) \int e^{-xt} \left[2n\left((\sinh t)^{2n-1} + (\sinh t)^{2n+1}\right) \right. \\ &\left. + (\sinh t)^{2n+1}\right] dt. \end{aligned} \tag{7}$$

Taking the limits in equation (7) from 0 to $\infty$, we obtain

$$\int_0^{\infty} e^{-xt} (\sinh t)^{2n+1} dt = \frac{-(2n+1)(2n)}{(2n+1)^2 - x^2} \int_0^{\infty} e^{-xt} (\sinh t)^{2n-1} dt. \tag{8}$$

Applying equation (8) $(n-1)$ times to the integral $\int_0^\infty e^{-xt}(\sinh t)^{2n-1}dt$, we obtain:

$$\begin{aligned}\int_0^\infty e^{-xt}(\sinh t)^{2n+1}dt &= \frac{(-1)^a(2n+1)(2n)(2n-1)(2n-2)\cdots.(3)(2)}{\left((2n+1)^2-x^2\right)\left((2n-1)^2-x^2\right)\cdots\left(3^2-x^2\right)} \\ &\cdot \int_0^\infty e^{-xt}\sinh t dt.\end{aligned} \tag{9}$$

The integral $\int_0^\infty e^{-xt}\sinh t dt$ can be calculated easily using integration by parts twicely to obtain:

$$\int_0^\infty e^{-xt}\sinh t dt = \frac{-1}{1-x^2}. \tag{10}$$

Substituting the fact in equation (10) into equation (9), we obtain:

$$\begin{aligned}\int_0^\infty e^{-xt}(\sinh t)^{2n+1}dt \\ = \frac{(-1)^{n+1}(2n+1)!}{\left((2n+1)^2-x^2\right)\left((2n-1)^2-x^2\right)\cdots\left(3-x^2\right)\left(1-x^2\right)}.\end{aligned} \tag{11}$$

Now, we express the solution of equation (5) in another form.

Using the power trigonometric formula deduced using De Moivre's formula, Euler's formula and the binomial theorem, see [11].

$$(\sinh(t))^{2n+1} = \frac{1}{4^n}\sum_{k=0}^{n}(-1)^k\binom{2n+1}{k}\sinh\left[(2n+1-2k)t\right]. \tag{12}$$

Substituting equation (12) into equation (5), we obtain:

$$\begin{aligned}\int_0^\infty e^{-xt}(\sinh t)^{2n+1}dt = \int_0^\infty e^{-xt}\frac{1}{4^n}\sum_{k=0}^{n}(-1)^k\binom{2n+1}{k} \\ \cdot \sinh\left[(2n+1-2k)t\right]dt.\end{aligned} \tag{13}$$

Therefore, by changing the order of the integral and the sum in equation (13), we obtain:

$$\begin{aligned}\int_0^\infty e^{-xt}(\sinh t)^{2n+1}dt = \frac{1}{4^n}\sum_{k=0}^{n}(-1)^k\binom{2n+1}{k} \\ \cdot \int_0^\infty e^{-xt}\sinh\left[(2n+1-2k)t\right]dt.\end{aligned} \tag{14}$$

To evaluate the integral $\int_0^\infty e^{-xt}\sinh\left[(2n+1-2k)t\right]dt$, we apply integration by parts twicely to obtain:

$$\int_0^\infty e^{-xt}\sinh\left[(2n+1-2k)t\right]dt = \frac{(2n+1-2k)}{x^2-(2n+1-2k)^2}. \tag{15}$$

Substituting the result in equation (15) into equation (14), we obtain:

$$\begin{aligned}\int_0^\infty e^{-xt}(\sinh t)^{2n+1}dt = \frac{1}{4^n}\sum_{k=0}^{n}(-1)^k\binom{2n+1}{k} \\ \cdot \frac{2n+1-2k}{x^2-(2n+1-2k)^2}.\end{aligned} \tag{16}$$

Equating equation (16) with equation (11), this completes the proof of equation (4). □

**Lemma 2.** *Let $n \in \mathbb{N}$, then*

$$\begin{aligned}\frac{1}{x\left(2^2-x^2\right)\left(4^2-x^2\right)\cdots\left(4n^2-x^2\right)} \\ = \frac{1}{2^{2n}(2n)!}\left(\frac{1}{x}\binom{2n}{n} \right. \\ \left. + 2\sum_{k=0}^{n-1}\left((-1)^{k+n+1}\binom{2n}{k}\frac{x}{(2n-2k)^2-x^2}\right)\right)\end{aligned} \tag{17}$$

*Proof.* The proof is done by repeating the same process in proving lemma (1) but using the integral $\int_0^\infty e^{-xt}(\sinh t)^{2n}dt$, where $x > 0, x > 2n, n \in \mathbb{N}$. □

**Lemma 3.** *Let $n = 0, 1, \cdots,$ and $m = 1, 2, \cdots,$ then we have*

$$\begin{aligned}&\frac{1}{\left[\left(1^2-x^2\right)\left(3^2-x^2\right)\cdots\left((2n+1)^2-x^2\right)\right]\left[x\left(2^2-x^2\right)\left(4^2-x^2\right)\cdots\left(4m^2-x^2\right)\right]} \\ &= \frac{(-1)^{n+1}}{x\,2^{2m+2n}(m!)^2(2n+1)!}\sum_{s=0}^{n}(-1)^s\binom{2n+1}{s}\frac{2n+1-2s}{x^2-(2n+1-2s)^2} \\ &+ \frac{(-1)^n}{2^{2m+2n-1}(m!)^2(2n+1)!}\sum_{k=0}^{m-1}\sum_{s=0}^{n}(-1)^{k+m+s+1}\binom{2m}{k}\binom{2n+1}{s}\frac{x(2n+1-2s)}{\left(\left(x^2-(2m-2k)^2\right)\left(x^2-(2n+1-2s)^2\right)\right)}\end{aligned} \tag{18}$$

*Proof.* This is a direct result obtained by multiplying equation (4) by equation (17). ☐

**Lemma 4.** *The following formulas of improper integrals are created using Lemma (1–3):*

$$PV\int_0^\infty \frac{\cos(\theta x)}{(1-x^2)(9-x^2)\cdots\left((2n+1)^2-x^2\right)}dx = \frac{(-1)^{n+1}\pi}{2^{2n+1}(2n+1)!}\sum_{k=0}^{n}(-1)^k\binom{2n+1}{k}\sin(\theta(2k-1-2n)), \tag{19}$$

*for* $\theta \geq 0, n = 0, 1, \cdots$.

*Proof.* The formula is obtained by multiplying both sides of equation (4) by $\cos(\theta x)$, then taking the Cauchy principal value of integral for both sides from 0 to $\infty$, and using the well-known fact:
$PV\int_0^\infty(\cos(\theta x)/a^2-x^2)dx=(\pi/2a)\sin(a\theta)$.
where $a, \theta > 0$,

$$PV\int_0^\infty \frac{x\sin(\theta x)}{(1-x^2)(9-x^2)\cdots\left((2n+1)^2-x^2\right)}dx = \frac{(-1)^{n+1}\pi}{2^{2n+1}(2n+1)!}\sum_{k=0}^{n}(-1)^k\binom{2n+1}{k}\cdot(2k-1-2n)\cos(\theta(2k-1-2n)), \tag{20}$$

for $\theta > 0, n = 0, 1, \cdots$. ☐

*Proof.* The formula is obtained by differentiating both sides of equation (19) with respect to $\theta$.

$$PV\int_0^\infty \frac{\sin(\theta x)}{x(4-x^2)(16-x^2)\cdots(4n^2-x^2)}dx = \frac{\pi}{2^{2n+1}(2n)!}\left(\binom{2n}{n}+2\sum_{k=0}^{n-1}(-1)^{k+n}\binom{2n}{k}\cos(\theta(2n-2k))\right), \tag{21}$$

for $\theta > 0, n = 1, 2, \cdots$. ☐

*Proof.* The formula is obtained by multiplying both sides of equation (17) by $\sin(\theta x)$, then taking the Cauchy principal value of integral from 0 to $\infty$, and using the well-known fact:
$PV\int_0^\infty(\sin(\theta x)/x(a^2-x^2))dx=(\pi/2a^2)(1-\cos(a\theta))$,
where $\theta, a > 0$

$$PV\int_0^\infty \frac{\cos(\theta x)}{(4-x^2)(16-x^2)\cdots(4n^2-x^2)}dx = \frac{\pi}{2^{2n}(2n)!}\sum_{k=0}^{n-1}(-1)^{k+n+1}\binom{2n}{k}(2n-2k)\sin(\theta(2k-2n)), \tag{22}$$

for $\theta \geq 0, n = 1, 2, \cdots$. ☐

*Proof.* The formula is obtained by differentiating both sides of equation (21) with respect to $\theta$. ☐

**Lemma 5.** *Let* $\theta > 0$ *and* $n = 0, 1, \cdots$, $m = 1, 2, \cdots$. *Then, we have the following improper integrals:*

$$\begin{aligned} PV\int_0^\infty &\frac{\sin(\theta x)}{\left((1-x^2)(9-x^2)\cdots\left((2n+1)^2-x^2\right)\right)\left(x(4-x^2)(16-x^2)\cdots(4m^2-x^2)\right)}dx \\ &= \frac{(-1)^n\pi}{2^{2m+2n+1}(m!)^2(2n+1)!}\sum_{s=0}^{n}(-1)^s\binom{2n+1}{s}\frac{(1-\cos(\theta(2n+1-2s)))}{(2n+1-2s)} \\ &+ \frac{(-1)^n\pi}{2^{2m+2n}(m!)^2(2n+1)!}\sum_{k=0}^{m-1}\sum_{s=0}^{n}(-1)^{k+m+s+1}\binom{2m}{k}\binom{2n+1}{s} \\ &\cdot\frac{(2n+1-2s)(\cos(\theta(2m-2k))-\cos(\theta(2n+1-2s)))}{(2m-2k)^2-(2n+1-2s)^2}. \end{aligned} \tag{23}$$

*Proof.* The formula is obtained by multiplying both sides of equation (18) by sin $(\theta x)$, then taking the Cauchy principal value of integral, and using the well-known facts:

$$PV\int_0^\infty \frac{\sin(\theta x)}{x(a^2-x^2)}dx = \frac{\pi}{2a^2}(1-\cos(a\theta)),$$

$$PV\int_0^\infty \frac{x\sin(\theta x)}{a^2-x^2}dx = -\frac{\pi}{2}\cos(a\theta), \quad (24)$$

where $\theta, a > 0$.

$$\begin{aligned}&\int_0^\infty \frac{\cos(\theta x)}{\left((1-x^2)(9-x^2)\cdots\left((2n+1)^2-x^2\right)\right)\left(x(4-x^2)(16-x^2)\cdots(4m^2-x^2)\right)}dx\\ &= \frac{(-1)^n\pi}{2^{2m+2n+1}(m!)^2(2n+1)!}\sum_{s=0}^{n}(-1)^s\binom{2n+1}{s}\sin(\theta(2n+1-2s)) + \frac{(-1)^n\pi}{2^{2m+2n}(m!)^2(2n+1)!}\sum_{k=0}^{m-1}\sum_{s=0}^{n}(-1)^{k+m+s+1}\\ &\quad\cdot\binom{2m}{k}\binom{2n+1}{s}(2n+1-2s)\frac{(2n+1-2s)\sin(\theta(2n+1-2s))-(2m-2k)\sin(\theta(2m-2k))}{(2m-2k)^2-(2n+1-2s)^2}.\end{aligned} \quad (25)$$

□

*Proof.* The formula is obtained by differentiating both sides of equation (23) with respect to $\theta$. □

## 3. New General Theorems

In this section, we present new master theorems to help mathematicians, engineers, and physicists solve complicated improper integrals. To obtain our goal, we present some facts about analytic functions [8, 10, 13–15].

Assuming that $f$ is an analytic function in a disc $D$ centered at $\alpha$, then using Taylor's expansion, where $\alpha$, $\beta$, and $\theta$ are real constants, we have

$$f(z) = \sum_{k=0}^{\infty}\frac{f^{(k)}(\alpha)}{k!}(z-\alpha)^k. \quad (26)$$

Substituting $z = \alpha + \beta e^{i\theta x}$ into $f(z)$, where $\beta$ is not completely arbitrary, since it must be smaller than the radius of $D$, we obtain

$$f\left(\alpha+\beta e^{i\theta x}\right) = \sum_{k=0}^{\infty}\frac{f^{(k)}(\alpha)}{k!}\beta^k e^{i\theta kx}. \quad (27)$$

Using the formulas

$$e^{i\theta x} + e^{-i\theta x} = 2\cos(\theta x),\ e^{i\theta x} - e^{-i\theta x} = 2i\sin(\theta x), \quad (28)$$

one can obtain

$$\begin{aligned}&\frac{1}{2}\left(f\left(\alpha+\beta e^{i\theta x}\right)+f\left(\alpha+\beta e^{-i\theta x}\right)\right)\\ &= \frac{1}{2}\sum_{k=0}^{\infty}\frac{f^{(k)}(\alpha)}{k!}\beta^k\left(e^{i\theta kx}+e^{-i\theta kx}\right) = \sum_{k=0}^{\infty}\frac{f^{(k)}(\alpha)}{k!}\beta^k\cos(k\theta x)\\ &= f(\alpha)+f'(\alpha)\beta\cos(\theta x)+\frac{f''(\alpha)}{2!}\beta^2\cos(2\theta x)+\cdots.\end{aligned} \quad (29)$$

Similarly,

$$\begin{aligned}&\frac{1}{2i}\left(f\left(\alpha+\beta e^{i\theta x}\right)-f\left(\alpha+\beta e^{-i\theta x}\right)\right)\\ &= \frac{1}{2}\sum_{k=0}^{\infty}\frac{f^{(k)}(\alpha)}{k!}\beta^k\left(e^{i\theta kx}-e^{-i\theta kx}\right)\\ &= f'(\alpha)\beta\sin(\theta x)+\frac{f''(\alpha)}{2!}\beta^2\sin(2\theta x)+\cdots\\ &= \sum_{k=1}^{\infty}\frac{f^{(k)}(\alpha)}{k!}\beta^k\sin(k\theta x).\end{aligned} \quad (30)$$

The parameters in equations (29) and (30) can be modified in the following lemma.

**Lemma 6.** *Assume that $g(\alpha+z)$ is an analytic function that has the following series expansion:*

$$g(\alpha+z) = \sum_{k=0}^{\infty}M_k e^{-kz}, \quad (31)$$

*whether $z$ be real or imaginary, and $\sum_{k=0}^{\infty}M_k$ is absolutely convergent. Then,*

$$\begin{aligned}\frac{1}{2}(g(\alpha-i\theta x)+g(\alpha+i\theta x)) &= \frac{1}{2}\sum_{k=0}^{\infty}M_k\left(e^{ik\theta x}+e^{-ik\theta x}\right)\\ &= \sum_{k=0}^{\infty}M_k\cos(k\theta x),\end{aligned} \quad (32)$$

*and,*

$$\frac{1}{2i}(g(\alpha-i\theta x)-g(\alpha+i\theta x))=\frac{1}{2i}\sum_{k=1}^{\infty}M_k\left(e^{ik\theta x}-e^{-ik\theta x}\right)=\sum_{k=1}^{\infty}M_k\sin(k\theta x), \tag{33}$$

*where >0, $\alpha\in\mathbb{R}$, and x is any real number.*

The next part of this section includes the new master theorems that we establish. Moreover, we mention here that Cauchy's results in [3] are identical to our results with special choices of the parameters, as will be discussed later.

**Theorem 1.** *Let f be an analytic function in a disc D centered at α, where $\alpha\in\mathbb{R}$. Then, we have the following improper integral formula:*

$$PV\int_0^{\infty}\frac{f\left(\alpha+\beta e^{i\theta x}\right)+f\left(\alpha+\beta e^{-i\theta x}\right)}{(1-x^2)(9-x^2)\cdots\left((2n+1)^2-x^2\right)}dx=\frac{(-1)^{n+1}\pi}{i\,2^{2n+1}(2n+1)!}\sum_{s=0}^{n}(-1)^s\binom{2n+1}{s}\cdot\left(f\left(\alpha+\beta e^{i\theta(2s-1-2n)}\right)-f\left(\alpha+\beta e^{-i\theta(2s-1-2n)}\right)\right), \tag{34}$$

*where $\theta\ge 0, n=0,1,2,\cdots$.*

*Proof.* of let

$$I=PV\int_0^{\infty}\frac{f\left(\alpha+\beta e^{i\theta x}\right)+f\left(\alpha+\beta e^{-i\theta x}\right)}{(1-x^2)(9-x^2)\cdots\left((2n+1)^2-x^2\right)}dx. \tag{35}$$

Now, since $f$ is an analytic function around $\alpha$ and substituting the fact in equation (29) into equation (35), we obtain:

$$I=PV\int_0^{\infty}\frac{2\sum_{k=0}^{\infty}\left(f^{(k)}(\alpha)\beta^k/k!\right)\cos(\theta kx)}{(1-x^2)(9-x^2)\cdots\left((2n+1)^2-x^2\right)}dx. \tag{36}$$

Using Fubini's theorem, the integral yields a finite answer when the integral is replaced by its absolute value, i.e., converges in the Riemann sense. Thus, we can interchange the order of the integration and the summation to obtain:

$$I=2\sum_{k=0}^{\infty}\frac{f^{(k)}(\alpha)\beta^k}{k!}PV\int_0^{\infty}\frac{\cos(\theta kx)}{(1-x^2)(9-x^2)\cdots\left((2n+1)^2-x^2\right)}dx. \tag{37}$$

Substituting the fact in equation (19) into equation (37), we obtain:

$$I=2\sum_{k=0}^{\infty}\frac{f^{(k)}(\alpha)\beta^k}{k!}\frac{(-1)^{n+1}\pi}{2^{2n+1}(2n+1)!}\sum_{s=0}^{n}(-1)^s\cdot\binom{2n+1}{s}\sin(\theta k(2s-1-2n)). \tag{38}$$

Rewriting $\sin(\theta k(2s-1-2n))$ in the exponential form, equation (38) becomes

$$I=\sum_{k=0}^{\infty}\frac{f^{(k)}(\alpha)\beta^k}{k!}\frac{(-1)^{n+1}\pi}{2^{2n}(2n+1)!}\frac{1}{2i}\sum_{s=0}^{n}(-1)^s\binom{2n+1}{s}\cdot\left(e^{i\theta k(2s-1-2n)}-e^{-i\theta k(2s-1-2n)}\right). \tag{39}$$

The fact in equation (27) implies that equation (39) becomes

$$I=\frac{(-1)^{n+1}\pi}{i2^{2n+1}(2n+1)!}\sum_{s=0}^{n}(-1)^s\binom{2n+1}{s}\cdot\left(f\left(\alpha+\beta e^{i\theta(2s-1-2n)}\right)-f\left(\alpha+\beta e^{-i\theta(2s-1-2n)}\right)\right). \tag{40}$$

Hence, this completes the proof. □

**Theorem 2.** *Let f be an analytic function in a disc D centered at α, where $\alpha\in\mathbb{R}$. Then, we have the following improper integral formula:*

$$PV\int_0^{\infty}\frac{x\left(f\left(\alpha+\beta e^{i\theta x}\right)-f\left(\alpha+\beta e^{-i\theta x}\right)\right)}{i(1-x^2)(9-x^2)\cdots\left((2n+1)^2-x^2\right)}dx=\frac{(-1)^{n+1}\pi}{2^{2n+1}(2n+1)!}\sum_{s=0}^{n}(-1)^s\binom{2n+1}{s}\cdot(2s-1-2n)(\psi(s)+\phi(s)-\eta), \tag{41}$$

*where $\theta>0$, $n=0,1,2,\cdots$, $\psi(s)=f(\alpha+\beta e^{i\theta(2s-1-2n)})$, $\phi(s)=f(\alpha+\beta e^{-i\theta(2s-1-2n)})$, and $\eta=2f(\alpha)$.*

*Proof.* Let

$$I=PV\int_0^{\infty}\frac{x\left(f\left(\alpha+\beta e^{i\theta x}\right)-f\left(\alpha+\beta e^{-i\theta x}\right)\right)}{i(1-x^2)(9-x^2)\cdots\left((2n+1)^2-x^2\right)}dx. \tag{42}$$

Now, since $f$ is an analytic function around $\alpha$ and thus substituting the fact in equation (30) into equation (7), we obtain

$$I=2\sum_{k=1}^{\infty}\frac{f^{(k)}(\alpha)\beta^k}{k!}PV\int_0^{\infty}\frac{x\sin(\theta kx)}{(1-x^2)(9-x^2)\cdots\left((2n+1)^2-x^2\right)}dx. \tag{43}$$

Substituting the fact in equation (20) into equation (43), we obtain

$$I=\sum_{k=1}^{\infty}\frac{f^{(k)}(\alpha)\beta^k}{k!}\frac{(-1)^{n+1}\pi}{2^{2n}(2n+1)!}\sum_{s=0}^{n}(-1)^k\binom{2n+1}{k}\cdot(2k-1-2n)\cos(\theta(2k-1-2n)). \tag{44}$$

Rewriting $\cos(\theta(2k-1-2n))$ in an exponential form, equation (44) becomes

$$I=\sum_{k=1}^{\infty}\frac{f^{(k)}(\alpha)\beta^k}{k!}\frac{(-1)^{n+1}\pi}{2^{2n}(2n+1)!}\frac{1}{2}\sum_{s=0}^{n}(-1)^s\binom{2n+1}{s}\cdot(2s-1-2n)\left(e^{i\theta k(2s-1-2n)}+e^{-i\theta k(2s-1-2n)}\right). \tag{45}$$

The fact in equation (27) implies that equation (45) becomes

$$I=\frac{(-1)^{n+1}\pi}{2^{2n+1}(2n+1)!}\sum_{s=0}^{n}(-1)^s\binom{2n+1}{s}\cdot(2s-1-2n)(\psi(s)+\phi(s)-\eta), \tag{46}$$

where $\psi(s)=f(\alpha+\beta e^{i\theta(2s-1-2n)})$, $\phi(s)=f(\alpha+\beta e^{-i\theta(2s-1-2n)})$, and $\eta=2f(\alpha)$.

Hence, this completes the proof. □

**Theorem 3.** *Let f be an analytic function in a disc D centered at α, where α ∈ ℝ. Then, we have the following improper integral formula:*

$$PV\int_0^{\infty}\frac{f(\alpha+\beta e^{i\theta x})-f(\alpha+\beta e^{-i\theta x})}{ix(4-x^2)(16-x^2)\cdots(4n^2-x^2)}dx=\frac{\pi}{2^{2n}(2n)!}\left(\binom{2n}{n}(\varphi-\eta)+\sum_{s=0}^{n-1}(-1)^{s+n}\binom{2n}{s}(\psi(s)+\phi(s)-2\eta)\right), \tag{47}$$

*where $\theta\geq 0$, $n=1,2,\cdots$, $\psi(s)=f(\alpha+\beta e^{i\theta(2n-2s)})$, $\phi(s)=f(\alpha+\beta e^{-i\theta(2n-2s)})$, $\eta=f(\alpha)$, and $\varphi=f(\alpha+\beta)$.*

*Proof.* The proof of Theorem 3 can be obtained by similar arguments to Theorem 2 and using the fact (21) in Lemma 4. □

**Theorem 4.** *Let f be an analytic function in a disc D centered at α, where α ∈ ℝ. Then, we have the following improper integral formula:*

$$PV\int_0^{\infty}\frac{f(\alpha+\beta e^{i\theta x})+f(\alpha+\beta e^{-i\theta x})}{(4-x^2)(16-x^2)\cdots(4n^2-x^2)}dx=\frac{\pi}{2^{2n}i(2n)!}\sum_{s=0}^{n-1}(-1)^{s+n+1}\cdot\binom{2n}{s}(2n-2s)\left(f\left(\alpha+\beta e^{i\theta(2n-2s)}\right)-f\left(\alpha+\beta e^{-i\theta(2n-2s)}\right)\right), \tag{48}$$

*where $\theta\geq 0$ $n=1,2,\cdots$.*

*Proof.* The proof of Theorem 4 can be obtained by similar arguments to Theorem 2 and using the fact (22) in Lemma 4. □

**Theorem 5.** *Let f be an analytic function in a disc D centered at α, where α ∈ ℝ. Then, we have the following improper integral formula:*

$$\begin{aligned}PV\int_0^{\infty}&\frac{f(\alpha+\beta e^{i\theta x})-f(\alpha+\beta e^{-i\theta x})}{i\left((1-x^2)(9-x^2)\cdots\left((2n+1)^2-x^2\right)\right)(x(4-x^2)(16-x^2)\cdots(4m^2-x^2))}dx\\&=\frac{(-1)^n\pi}{2^{2m+2n}(m!)^2(2n+1)!}\sum_{s=0}^{n}(-1)^s\binom{2n+1}{s}\frac{(\varphi-1/2(\psi(s)+\phi(s)))}{(2n+1-2s)}\\&+\frac{(-1)^n\pi}{2^{2m+2n-1}(m!)^2(2n+1)!}\sum_{k=0}^{m-1}\sum_{s=0}^{n}(-1)^{k+m+s+1}\binom{2m}{k}\binom{2n+1}{s}(2n+1-2s)\frac{1/2((\gamma(k)+\lambda(k))-1/2(\psi(s)+\phi(s)))}{\left((2m-2k)^2-(2n+1-2s)^2\right)},\end{aligned} \tag{49}$$

*where* $\theta > 0$, $n = 0, 1, 2, \cdots, m = 1, 2, \cdots$, $\psi(s) = f(\alpha + \beta e^{i\theta(2n+1-2s)})$, $\phi(s) = f(\alpha + \beta e^{-i\theta(2n+1-2s)})$, $\gamma(k) = f(\alpha + \beta e^{i\theta(2m-2k)})$, $\lambda(k) = f(\alpha + \beta e^{-i\theta(2m-2k)})$, *and* $\varphi = f(\alpha + \beta)$.

*Proof.* Let

$$I = PV\int_0^\infty \frac{f\left(\alpha + \beta e^{i\theta x}\right) - f\left(\alpha + \beta e^{-i\theta x}\right)}{i\left((1-x^2)(9-x^2)\cdots\left((2n+1)^2 - x^2\right)\right)\left(x(4-x^2)(16-x^2)\cdots(4m^2-x^2)\right)}dx. \quad (50)$$

Now, since $f$ is an analytic function around $\alpha$ and substituting the fact in equation (30) into equation (50), we obtain

$$I = 2\sum_{j=1}^{\infty}\frac{f^{(j)}(\alpha)\beta^j}{j!} \quad (51)$$

$$PV\int_0^\infty \frac{\sin(\theta jx)}{i\left((1-x^2)(9-x^2)\cdots\left((2n+1)^2 - x^2\right)\right)\left(x(4-x^2)(16-x^2)\cdots(4m^2-x^2)\right)}dx. \quad (52)$$

Substituting the fact in equation (23) into equation (51), we obtain

$$I = 2\sum_{j=1}^{\infty}\frac{f^{(j)}(\alpha)\beta^j}{j!}(A+B), \quad (53)$$

where

$$A = \frac{(-1)^n\pi}{2^{2m+2n+1}(m!)^2(2n+1)!}\sum_{s=0}^{n}(-1)^s\binom{2n+1}{s} \cdot \frac{1-\cos(\theta j(2n+1-2s))}{2n+1-2s},$$

$$B = \frac{(-1)^n\pi}{2^{2m+2n}(m!)^2(2n+1)!}\sum_{k=0}^{m-1}\sum_{s=0}^{n}(-1)^{k+m+s+1} \cdot \binom{2m}{k}\binom{2n+1}{s}(2n+1-2s) \cdot \frac{\cos(\theta j(2m-2k)) - \cos(\theta j(2n+1-2s))}{(2m-2k)^2 - (2n+1-2s)^2}. \quad (54)$$

Rewriting $\cos(\theta j(2n+1-2s))$, and $\cos(\theta j(2m-2k))$ in the exponential forms and using equation (5), equation (53) becomes

$$\begin{aligned} &PV\int_0^\infty \frac{f\left(\alpha + \beta e^{i\theta x}\right) - f\left(\alpha + \beta e^{-i\theta x}\right)}{i\left((1-x^2)(9-x^2)\cdots\left((2n+1)^2 - x^2\right)\right)\left(x(4-x^2)(16-x^2)\cdots(4m^2-x^2)\right)}dx \\ &= \frac{(-1)^n\pi}{2^{2m+2n}(m!)^2(2n+1)!}\sum_{s=0}^{n}(-1)^s\binom{2n+1}{s}\frac{(\varphi - 1/2(\psi(s)+\phi(s)))}{(2n+1-2s)} \\ &+ \frac{(-1)^n\pi}{2^{2m+2n-1}(m!)^2(2n+1)!}\sum_{k=0}^{m-1}\sum_{s=0}^{n}(-1)^{k+m+s+1}\binom{2m}{k}\binom{2n+1}{s}(2n+1-2s)\frac{1/2((\gamma(k)+\lambda(k)) - 1/2(\psi(s)+\phi(s)))}{\left((2m-2k)^2 - (2n+1-2s)^2\right)}, \end{aligned} \quad (55)$$

where $\psi(s) = f(\alpha + \beta e^{i\theta(2n+1-2s)})$, $\phi(s) = f(\alpha + \beta e^{-i\theta(2n+1-2s)})$, $\gamma(k) = f(\alpha + \beta e^{i\theta(2m-2k)})$, $\lambda(k) = f(\alpha + \beta e^{-i\theta(2m-2k)})$, and $\varphi = f(\alpha + \beta)$.

Hence, this completes the proof of Theorem 5. □

**Theorem 6.** *Let $f$ be an analytic function in a disc $D$ centered at $\alpha$, where $\alpha \in \mathbb{R}$. Then, we have the following improper integral formula:*

$$PV\int_0^\infty \frac{f\left(\alpha+\beta e^{i\theta x}\right)+f\left(\alpha+\beta e^{-i\theta x}\right)}{\left((1-x^2)(9-x^2)\cdots\left((2n+1)^2-x^2\right)\right)\left((4-x^2)(16-x^2)\cdots(4m^2-x^2)\right)}dx$$
$$= \frac{(-1)^n\pi}{2^{2m+2n}(m!)^2(2n+1)!}\sum_{s=0}^{n}(-1)^s\binom{2n+1}{s}\frac{1}{2i}(\psi(s)-\phi(s)) + \frac{(-1)^n\pi}{2^{2m+2n-1}(m!)^2(2n+1)!}\sum_{k=0}^{m-1}\sum_{s=0}^{n}(-1)^{k+m+s+1}\binom{2m}{k}\binom{2n+1}{s}$$
$$\cdot (2n+1-2s)\,\frac{(2n+1-2s)/2i(\psi(s)-\phi(s))-(2m-2k)/2i(\gamma(k)-\lambda(k))}{\left((2m-2k)^2-(2n+1-2s)^2\right)}, \tag{56}$$

*where* $\geq 0$, $n = 0, 1, 2, \cdots$, $m = 1, 2, \cdots$, $\psi(s) = f(\alpha + \beta e^{i\theta(2n+1-2s)})$, $\phi(s) = f(\alpha + \beta e^{-i\theta(2n+1-2s)})$, $\gamma(k) = f(\alpha + \beta e^{i\theta(2m-2k)})$, *and* $\lambda(k) = f(\alpha + \beta e^{-i\theta(2m-2k)})$.

*Proof.* The proof of Theorem 6 can be obtained by similar arguments to Theorem 5 and using the fact (25) in Lemma 5. □

The following table illustrates some corollaries of the previous theorems with special cases under the assumption in Lemma 6. We introduce the principal value of some improper integrals, which are special cases of the proposed theorems.

## 4. Numerical Applications

In this section, we present the results, applications, and observations of the proposed theorems. We also show that the simple cases of the master theorems are identical to the results obtained by Cauchy in his memoirs using Residue Theorem 4. Also, some examples on difficult integrals that cannot be treated directly by usual methods. In this section, we show the applicability of our results in handling such problems.

*4.1. Some Remarks on the Theorems*

*Remark 1.* Letting $\alpha = 0$ and $n = 1$ in Theorem 3, we obtain

$$PV\int_0^\infty \frac{f\left(\beta e^{i\theta x}\right)-f\left(\beta e^{-i\theta x}\right)}{ix(4-x^2)}dx = \frac{\pi}{4}\left(f(\beta)-\frac{1}{2}\left(f\left(\beta e^{2i\theta}\right)+f\left(\beta e^{-2i\theta}\right)\right)\right), \tag{57}$$

where $\theta > 0$.

Letting $x/2 = y$,

$$\frac{1}{4}PV\int_0^\infty \frac{f\left(\beta e^{2i\theta y}\right)-f\left(\beta e^{-2i\theta y}\right)}{iy(1-y^2)}dy = \frac{\pi}{4}\left(f(\beta)-\frac{1}{2}\left(f\left(\beta e^{2i\theta}\right)+f\left(\beta e^{-2i\theta}\right)\right)\right). \tag{58}$$

Letting $2\theta = \varphi$,

$$PV\int_0^\infty \frac{f\left(\beta e^{i\varphi y}\right)-f\left(\beta e^{-i\varphi y}\right)}{iy(1-y^2)}dy = \pi\left(f(\beta)-\frac{1}{2}\left(f\left(\beta e^{i\varphi}\right)+f\left(\beta e^{-i\varphi}\right)\right)\right). \tag{59}$$

This result does not appear in [4–6, 11].

The following table presents comparisons to Cauchy's results, which illustrate the relation between our theorems and the results obtained by Cauchy; that is, some of Cauchy's results become simple cases of our general theorems.

*4.2. Generating Improper Integrals.* In this section, we show the mechanism of generating an infinite number of integrals by choosing the function $f(z)$ and finding the real or imaginary part. It is worth noting that some of these integrals with special cases appear in [31–34] when solving some applications related to finding Green's function, one-dimensional vibrating string problems, wave motion in elastic solids, and using Fourier cosine and Fourier Sine transforms.

To illustrate the idea, we show some general examples that are applied on Theorems 1, 2, and 3, as follows:

(1) Setting $f(z)=z^m$, $m\in\mathbb{R}^+$:

(i) Using Theorem 1 and setting $\alpha=0$, $\beta=1$ we have:

$$f\left(e^{i\theta x}\right)+f\left(e^{-i\theta x}\right)=e^{i\theta mx}+e^{-i\theta mx}=2\cos(\theta mx). \tag{60}$$

Thus,

$$\int_0^\infty \frac{2\cos(\theta x)}{(1-x^2)(9-x^2)\cdots\left((2n+1)^2-x^2\right)}dx = \frac{(-1)^{n+1}\pi}{2^{2n}(2n+1)!}\sum_{k=0}^{n}(-1)^k\binom{2n+1}{k} \cdot \sin(m\theta(2k-1-2n)), \tag{61}$$

for $\theta\geq 0$, $n=0,1,\cdots$.

Setting $m=1$, the obtained integral is a Fourier cosine transform [31, 32] of the function $f(t)=1/((1-t^2)(9-t^2)\cdots((2n+1)^2-t^2))$

(ii) Using Theorem 3 and setting $\alpha=0$, $\beta=1$ we have:

$$\frac{1}{i}\left(f\left(e^{i\theta x}\right)-f\left(e^{-i\theta x}\right)\right)=\frac{1}{i}\left(e^{i\theta mx}+e^{-i\theta mx}\right) = 2\sin(\theta mx). \tag{62}$$

Thus,

$$PV\int_0^\infty \frac{2\sin(\theta mx)}{x(4-\mathrm{x}^2)(16-x^2)\cdots(4n^2-x^2)}dx = \frac{\pi}{2^{2n}(2n)!}\left(\binom{2n}{n}+2\sum_{k=0}^{n-1}(-1)^{k+n} \cdot \binom{2n}{k}\cos(\theta m(2n-2k))\right), \tag{63}$$

for $\theta>0$, $n=1,2,\cdots$.

Setting $m=1$, the obtained integral is a Fourier sine transform [31, 32] of the function $f(t)=1/(t(4-t^2)(16-t^2)\cdots(4n^2-t^2))$.

Setting $f(z)=e^z$

(i) Using Theorem 1, we have:

$$f\left(\alpha+\beta e^{i\theta x}\right)+f\left(\alpha+\beta e^{-i\theta x}\right) = e^{\alpha+\beta e^{i\theta x}}+e^{\alpha+\beta e^{-i\theta x}}=2e^{\alpha+\beta\cos(\theta x)}\cos(\beta\sin(\theta x))$$

$$\begin{aligned}PV\int_0^\infty &\frac{2e^{\alpha+\beta\cos(\theta x)}\cos(\beta sin(\theta x))}{(1-x^2)(9-x^2)\cdots\left((2n+1)^2-x^2\right)}dx \\ &= \frac{(-1)^{n+1}\pi}{2^{2n+1}(2n+1)!}\sum_{s=0}^{n}(-1)^s\binom{2n+1}{s} \cdot\left(\frac{e^\alpha}{i}\left(e^{\beta e^{i\theta(2s-1-2n)}}-e^{\beta e^{-i\theta(2s-1-2n)}}\right)\right) \\ &= \frac{(-1)^{n+1}\pi}{2^{2n}(2n+1)!}\sum_{s=0}^{n}(-1)^s\binom{2n+1}{s} \cdot\left(e^\alpha(\sin(\beta\sin(\theta(2s-1-2n)))\right. \\ &\cdot(\sinh(\beta\cos(\theta(2s-1-2n))) \\ &\left.+\cosh(\beta\cos(\theta(2s-1-2n)))))\right),\end{aligned} \tag{64}$$

where $\theta\geq 0$ and $n=0,1,2,\cdots$

(ii) Using Theorem 3, we have:

$$\begin{aligned}\frac{1}{i}&\left(f\left(\alpha+\beta e^{i\theta x}\right)-f\left(\alpha+\beta e^{-i\theta x}\right)\right) \\ &=\frac{1}{i}\left(e^{\alpha+\beta e^{i\theta x}}-e^{\alpha+\beta e^{-i\theta x}}\right) \\ &=2e^{\alpha+\beta\cos(\theta x)}\sin(\beta\sin(\theta x))\end{aligned} \tag{65}$$

Thus,

$$\begin{aligned}PV\int_0^\infty &\frac{x\left(2e^{\alpha+\beta\cos(\theta x)}\sin(\beta sin(\theta x))\right)}{(1-x^2)(9-x^2)\cdots\left((2n+1)^2-x^2\right)}dx \\ &= \frac{(-1)^{n+1}\pi}{2^{2n+1}(2n+1)!}\sum_{s=0}^{n}(-1)^s\binom{2n+1}{s} \\ &\cdot(2s-1-2n)\left(e^{\alpha+\beta e^{i\theta(2s-1-2n)}}+e^{\alpha+\beta e^{-i\theta(2s-1-2n)}}-2e^\alpha\right),\end{aligned} \tag{66}$$

where $\theta>0$, $n=0,1,2,\cdots$.

Setting $f(z)=\sin hz$

(iii) Using Theorem 1, we have:

$$\begin{aligned}f&\left(\alpha+\beta e^{i\theta x}\right)+f\left(\alpha+\beta e^{-i\theta x}\right) \\ &=\sinh\left(\alpha+\beta e^{i\theta x}\right)+\sinh\left(\alpha+\beta e^{-i\theta x}\right) \\ &=2\cos(\beta\sin(\theta x))\sinh(\alpha+\beta\cos(\theta x)).\end{aligned} \tag{67}$$

Thus,

$$PV\int_0^\infty \frac{2\cos(\beta\sin(\theta x))\sinh(\alpha+\beta\cos(\theta x))}{(1-x^2)(9-x^2)\cdots((2n+1)^2-x^2)}dx = \frac{(-1)^{n+1}\pi}{2^{2n+1}(2n+1)!}\sum_{s=0}^{n}(-1)^s\binom{2n+1}{s}\frac{1}{2i} \cdot\left(\sinh\left(\alpha+\beta e^{i\theta(2s-1-2n)}\right) - \sinh\left(\alpha+\beta e^{-i\theta(2s-1-2n)}\right)\right). \tag{68}$$

(iv) Using Theorem 2, we have:

$$\frac{1}{i}\left(f\left(\alpha+\beta e^{i\theta x}\right)-f\left(\alpha+\beta e^{-i\theta x}\right)\right) = \frac{1}{i}\left(\sinh\left(\alpha+\beta e^{i\theta x}\right)-\sinh\left(\alpha+\beta e^{-i\theta x}\right)\right) = 2\sin(\beta\sin(\theta x))\cosh(\alpha+\beta\cos(\theta x)). \tag{69}$$

Thus,

$$PV\int_0^\infty \frac{2\sin(\beta\sin(\theta x))\cosh(\alpha+\beta\cos(\theta x))}{x(4-x^2)(16-x^2)\cdots(4n^2-x^2)}dx = \frac{\pi}{2^{2n}(2n)!}\left(\binom{2n}{n}(\sinh(\alpha+\beta)-\sinh(\alpha))\right) + \frac{\pi}{2^{2n}(2n)!}\sum_{s=0}^{n-1}(-1)^{s+n}\binom{2n}{s} \cdot\left(\sinh\left(\alpha+\beta e^{i\theta(2n-2s)}\right) + \sinh\left(\alpha+\beta e^{-i\theta(2n-2s)}\right)-2\sinh(\alpha)\right) \tag{70}$$

where $\theta>0, n=1,2,\cdots$.

Setting $f(z)=\cos(e^z)$

(i) Using Theorem 1, we have:

$$f\left(\alpha+\beta e^{i\theta x}\right)+f\left(\alpha+\beta e^{-i\theta x}\right) = \cos\left(e^{\alpha+\beta e^{i\theta x}}\right)+\cos\left(e^{\alpha+\beta e^{-i\theta x}}\right) = 2\cos\left(e^{\alpha+\beta\cos(\theta x)}\cos(\beta sin(\theta x))\right) \cdot\cosh\left(\sin(\beta sin(\theta x))e^{\alpha+\beta\cos(\theta x)}\right) \tag{71}$$

Thus,

$$PV\int_0^\infty \frac{2\cos\left(e^{\alpha+\beta\cos(\theta x)}\cos(\beta sin(\theta x))\right)\cosh\left(\sin(\beta sin(\theta x))e^{\alpha+\beta\cos(\theta x)}\right)}{(1-x^2)(9-x^2)\cdots((2n+1)^2-x^2)}dx = \frac{(-1)^{n+1}\pi}{2^{2n}(2n+1)!}\sum_{s=0}^{n}(-1)^s\binom{2n+1}{s}\frac{1}{2i}\left(\cos\left(e^{\alpha+\beta e^{i\theta(2s-1-2n)}}\right)-\cos\left(e^{\alpha+\beta e^{-i\theta(2s-1-2n)}}\right)\right). \tag{72}$$

Setting $f(z)=\ln(1+z)$

(i) Using Theorem 1, we have:

$$f\left(1+\alpha+\beta e^{i\theta x}\right)+f\left(1+\alpha+\beta e^{-i\theta x}\right) = \ln\left(1+\alpha+\beta e^{i\theta x}\right)+\ln\left(1+\alpha+\beta e^{-i\theta x}\right) = \ln\left((\alpha+1)^2+\beta^2+2(\alpha+1)\beta\cos(\theta x)\right). \tag{73}$$

Thus,

$$PV\int_0^\infty \frac{\ln\left((\alpha+1)^2+\beta^2+2(\alpha+1)\beta\cos(\theta x)\right)}{(1-x^2)(9-x^2)\cdots((2n+1)^2-x^2)}dx = \frac{(-1)^{n+1}\pi}{2^{2n}(2n+1)!}\sum_{s=0}^{n}(-1)^s\binom{2n+1}{s} \cdot\frac{1}{2i}\left(\ln\left(1+\alpha+\beta e^{i\theta(2s-1-2n)}\right)-\ln\left(1+\alpha+\beta e^{-i\theta(2s-1-2n)}\right)\right) \tag{74}$$

(ii) Setting $\alpha = 0, \beta = 1$, we have:

$$f\left(e^{i\theta x}\right) + f\left(e^{-i\theta x}\right) = \ln\left(1 + e^{i\theta x}\right) + \ln\left(1 + e^{-i\theta x}\right) = 2\ln\left|2\cos\left(\frac{\theta x}{2}\right)\right|. \tag{75}$$

Thus,

$$PV\int_0^\infty \frac{2\ln|2\cos(\theta x/2)|}{(1-x^2)(9-x^2)\cdots\left((2n+1)^2-x^2\right)}dx = \frac{(-1)^{n+1}\pi}{2^{2n}(2n+1)!}\sum_{s=0}^{n}(-1)^s\binom{2n+1}{s} \cdot \frac{1}{2i}\left(\ln\left(1+e^{i\theta(2s-1-2n)}\right) - \ln\left(1+e^{-i\theta(2s-1-2n)}\right)\right) \tag{76}$$

*4.3. Solving Improper Integrals.* In this section, some applications on complicated problems are introduced and solved directly using a particular case of our new theorems. We note that the Mathematica and Maple software cannot solve such examples.

*Example 1.* Evaluate the following integral:

$$PV\int_0^\infty \frac{xe^{e^{\cos(\theta x)}\cos(\sin(\theta x))}\sin\left(\sin(\sin(\theta x))e^{\cos(\theta x)}\right)}{1-x^2}dx, \tag{77}$$

Solution Using Theorem 2, and let $f(z) = e^{e^z}$, we have

$$\frac{1}{i}\left(f\left(\alpha+\beta e^{i\theta x}\right) - f\left(\alpha+\beta e^{-i\theta x}\right)\right) = 2e^{e^{\alpha+\beta\cos(\theta x)}\cos(\beta\sin(\theta x))}\sin\left(\sin(\beta\sin(\theta x))e^{\alpha+\beta\cos(\theta x)}\right). \tag{78}$$

Thus,

$$PV\int_0^\infty \frac{xe^{e^{\alpha+\beta\cos(\theta x)}\cos(\beta\sin(\theta x))}\sin\left(\sin(\beta\sin(\theta x))e^{\alpha+\beta\cos(\theta x)}\right)}{(1-x^2)(9-x^2)\cdots\left((2n+1)^2-x^2\right)}dx = \frac{(-1)^{n+1}\pi}{2^{2n+1}(2n+1)!}\sum_{s=0}^{n}(-1)^s\binom{2n+1}{s}(2s-1-2n) \cdot \left(e^{e^{\alpha+\beta e^{i\theta(2s-1-2n)}}} + e^{e^{\alpha+\beta e^{-i\theta(2s-1-2n)}}} - 2e^{e^{\alpha}}\right). \tag{79}$$

Setting $n = 0, \alpha = 0, \beta = 1$,, we obtain

$$PV\int_0^\infty \frac{xe^{e^{\cos(\theta x)}\cos(\sin(\theta x))}\sin\left(\sin(\sin(\theta x))e^{\cos(\theta x)}\right)}{1-x^2}dx = \frac{\pi}{2}\left(e^{e^{e^{i\theta}}} + e^{e^{e^{-i\theta}}} - 2e\right). \tag{80}$$

which is simplified to $\pi(e^{e^{\cos(\theta)}\cos(\sin(\theta))}\cos(\sin(\sin(\theta))e^{\cos(\theta)}) - e)$

*Example 2.* Evaluate the following integral:

$$PV\int_0^\infty \frac{\left(e^{b\tan^{-1}(\theta x)} + e^{-b\tan^{-1}(\theta x)}\right)\cos\left((b/2)\ln\left(1+\theta^2x^2\right)\right)}{\left(1^2-x^2\right)\left(3^2-x^2\right)}dx, \tag{81}$$

where $\theta > 0, b \in \mathbb{R}$.

Solution: using result 1 in Table 1 and setting $\alpha = 0$, $n = 1$ and $g(z) = \cos(b\ln(1+z))$. We have

$$g(-i\theta x) + g(i\theta x) = \cos(b\ln(1-i\theta x)) + \cos(b\ln(1+i\theta x)) = \cos\left(\frac{b}{2}\ln\left(1+\theta^2x^2\right) + ib\tan^{-1}(\theta x)\right) + \cos\left(\frac{b}{2}\ln\left(1+\theta^2x^2\right) - ib\tan^{-1}(\theta x)\right) = 2\cosh\left(b\tan^{-1}(\theta x)\right)\cos\left(\frac{b}{2}\ln\left(1+\theta^2x^2\right)\right). \tag{82}$$

Therefore, we obtain

$$PV\int_0^\infty \frac{2\cosh\left(b\tan^{-1}(\theta x)\right)\cos\left((b/2)\ln\left(1+\theta^2x^2\right)\right)}{(1-x^2)(9-x^2)}dx = \frac{\pi}{i\,48}(g(3i\theta) - g(-3i\theta) - 3(g(i\theta) - g(-i\theta))) = \frac{\pi}{24}\left(3\sin\left(\frac{1}{2}b\ln\left(\theta^2x^2+1\right)\right)\sinh\left(b\tan^{-1}(\theta x)\right) - \sin\left(\frac{1}{2}b\ln\left(9\theta^2x^2+1\right)\right)\sinh\left(b\tan^{-1}(3\theta x)\right)\right). \tag{83}$$

Setting $b = 1, \theta = 1$, we obtain

$$PV\int_0^\infty \frac{2\cosh\left(\tan^{-1}(x)\right)\cos\left(\ln\left(\sqrt{1+x^2}\right)\right)}{(1-x^2)(9-x^2)}dx = \frac{\pi}{24}\left(3\sin\left(\ln\left(\sqrt{x^2+1}\right)\right)\sinh\left(\tan^{-1}(x)\right) - \sin\left(\ln\left(\sqrt{9x^2+1}\right)\right)\sinh\left(\tan^{-1}(3x)\right)\right).$$

Table 1: The Cauchy principal value of improper integrals that have the series representation in Equation (31).

| | $f(x)$ | $PV\int_0^{\infty} f(x)dx$ | Conditions | No. of theorem |
|---|---|---|---|---|
| 1. | $\frac{g(\alpha - i\theta x) + g(\alpha + i\theta x)}{(1-x^2)(9-x^2)\cdots\left((2n+1)^2 - x^2\right)}$ | $\frac{(-1)^{n+1}\pi}{i2^{2n+1}(2n+1)!}\sum_{s=0}^{n}(-1)^s\binom{2n+1}{s}\left(g(\alpha - i\theta(2s-1-2n)) - g(\alpha + i\theta(2s-1-2n))\right),$ | $\theta \geq 0,$ $n = 0, 1, 2, \cdots$ | Theorem 1 |
| 2. | $\frac{x(g(\alpha - i\theta x) - g(\alpha + i\theta x))}{i(1-x^2)(9-x^2)\cdots\left((2n+1)^2 - x^2\right)}$ | $\frac{(-1)^{n+1}\pi}{2^{2n+1}(2n+1)!}\sum_{s=0}^{n}(-1)^s\binom{2n+1}{s}(2s-1-2n)(\psi(s) + \phi(s) - \eta),$ Where $\psi(s) = g(\alpha - i\theta(2s-1-2n)), \phi(s) = g(\alpha + i\theta(2s-1-2n))$,and $\eta = 2g(\alpha)$. | $\theta > 0,$ $n = 0, 1, 2, \cdots$ | Theorem 2 |
| 3. | $\frac{g(\alpha - i\theta x) + g(\alpha + i\theta x)}{(4-x^2)(16-x^2)\cdots(4n^2 - x^2)}$ | $\frac{\pi}{2^{2n}i(2n)!}\sum_{s=0}^{n-1}(-1)^{s+n+1}\binom{2n}{s}(2n-2s)\left(g(\alpha - i\theta(2n-2s)) - g(\alpha + i\theta(2n-2s))\right).$ | $\theta \geq 0,$ $n = 1, 2, \cdots$ | Theorem 4 |
| 4. | $\frac{\frac{g(\alpha - i\theta x) + f(\alpha + i\theta x)}{\left((1-x^2)(9-x^2)\cdots\left((2n+1)^2 - x^2\right)\right)}}{\frac{1}{((4-x^2)(16-x^2)\cdots(4m^2 - x^2))}}$ | $\frac{(-1)^n\pi}{2^{2m+2n}(m!)^2(2n+1)!}\sum_{s=0}^{n}(-1)^s\binom{2n+1}{s}\frac{1}{2i}(\psi(s) - \phi(s)) + \frac{(-1)^n\pi}{2^{2m+2n-1}(m!)^2(2n+1)!}\sum_{k=0}^{m-1}\sum_{s=0}^{n}(-1)^{k+m+s+1}\binom{2m}{k}\binom{2n+1}{s}$ $\cdot(2n+1-2s)\,\frac{((2n+1-2s)/2i(\psi(s) - \phi(s)) - (2m-2k)/2i(\gamma(k) - \lambda(k)))}{\left((2m-2k)^2 - (2n+1-2s)^2\right)},$ Where $\psi(s) = g(\alpha - i\theta(2n+1-2s)), \phi(s) = f(\alpha + i\theta(2n+1-2s)), \gamma(k) = g(\alpha - i\theta(2m-2k))$, and $\lambda(k) = g$ $(\alpha + i\theta(2m-2k))$. | $\theta \geq 0,$ $n = 0, 1, 2, \cdots$ $m = 1, 2, \cdots$ | Theorem 6 |

Table 2: Remarks on principal values integrals, where $\theta > 0$.

| | Conditions | Theorem | $g(x)$ | $PV\int_0^{\infty} g(x)dx$ | Remarks |
|---|---|---|---|---|---|
| 1 | $\alpha = 0$, $\beta = 1$, and $n = 0$ | 1 | $\dfrac{f(e^{i\theta x}) + f(e^{-i\theta x})}{(1 - x^2)}$ | $\dfrac{\pi}{2i}\left(f\left(e^{i\theta}\right) - f\left(e^{-i\theta}\right)\right)$ | Which is identical to Cauchy's theorem, obtained [4] (P.62 formula (9)). |
| 2 | $\alpha = 0$, $\beta = 1$, and $n = 1$ | 2 | $\dfrac{x(f(e^{i\theta x}) - f(e^{-i\theta x}))}{i(1 - x^2)}$ | $\pi\left(f(0) - \dfrac{1}{2}\left(f\left(e^{i\theta}\right) + f\left(e^{-i\theta}\right)\right)\right)$ | Cauchy also made a mistake in this result, see [4] (P.62 formula (9)). He did not correct this mistake in the next memoir see [5, 6]. |
| 3 | $n = 1$ | 1 | $\dfrac{f(\alpha + \beta e^{i\theta x}) + f(\alpha + \beta e^{-i\theta x})}{(1^2 - x^2)(3^2 - x^2)}$ | $\dfrac{\pi}{i\,48}\left(f\left(\alpha + \beta e^{-3i\theta}\right) - f\left(\alpha + \beta e^{3i\theta}\right) - 3\left(f\left(\alpha + \beta e^{-i\theta}\right) - f\left(\alpha + \beta e^{i\theta}\right)\right)\right)$ | This result does not appear in [4–6, 11]. |
| 4 | $n = 1$ | 2 | $\dfrac{x(f(\alpha + \beta e^{i\theta x}) - f(\alpha + \beta e^{-i\theta x}))}{i(1^2 - x^2)(3^2 - x^2)}$ | $\dfrac{\pi}{16}\left(f(\alpha + \beta e^{-i\theta}) + f(\alpha + \beta e^{i\theta}) - f\left(\alpha + \beta e^{-3i\theta}\right) - f\left(\alpha + \beta e^{3i\theta}\right)\right)$ | ……… |
| 5 | $n = 2$ | 3 | $\dfrac{f(\alpha + \beta e^{i\theta x}) - f(\alpha + \beta e^{-i\theta x})}{ix(2^2 - x^2)(4^2 - x^2)}$ | $\dfrac{\pi}{384}\left(6(f(\alpha + \beta)) + f\left(\alpha + \beta e^{4i\theta}\right) + f\left(\alpha + \beta e^{-4i\theta}\right) - 4\left(f\left(\alpha + \beta e^{2i\theta}\right) + f\left(\alpha + \beta e^{-2i\theta}\right)\right)\right)$ | ……… |
| 6 | $n = 2$ | 4 | $\dfrac{f(\alpha + \beta e^{i\theta x}) + f(\alpha + \beta e^{-i\theta x})}{(2^2 - x^2)(4^2 - x^2)}$ | $\dfrac{\pi}{96\,i}\left(2\left(f(\alpha + \beta e^{2i\theta}) - f(\alpha + \beta e^{-2i\theta})\right) - \left(f(\alpha + \beta e^{4i\theta}) - f(\alpha + \beta e^{-4i\theta})\right)\right)$ | ……… |

$$\frac{\pi}{96\,i}\left(2\left(\tan^{-1}e^{2i\theta}\right)^2-\left(\tan^{-1}e^{-2i\theta}\right)^2-\left(\left(\tan^{-1}e^{2i\theta}\right)^2-\left(\tan^{-1}e^{-2i\theta}\right)^2\right)\right) \quad (84)$$

*Example 3.* Evaluate the following integral:

$$PV\int_0^{\infty}\frac{\ln^2|\tan((\theta x/2)-(\pi/4))|}{(2^2-x^2)(4^2-x^2)}dx, \quad (85)$$

where $\theta>0$.

Solution: using Theorem 4 and setting $\alpha=0$, $\beta=1$, $n=2$ or using Remark 6 Table 2 and setting $=0$, $\beta=1$, we set

$$f(z)=\left(\tan^{-1}z\right)^2=(-1/4)\ln^2((1-iz)/(1+iz)). \quad (86)$$

Therefore, we have $f(e^{i\theta x})=-1/4\ln^2((1-ie^{i\theta x})/(1+ie^{i\theta x}))$, and $f(e^{i\theta x})+f(e^{-i\theta x})=2\ \mathrm{Re}\,f(e^{i\theta x})$

Thus, we obtain

$$\begin{aligned}
&PV\int_0^{\infty}\frac{f(e^{i\theta x})+f(e^{-i\theta x})}{(2^2-x^2)(4^2-x^2)}dx\\
&=\frac{-1}{4}PV\int_0^{\infty}\frac{\ln^2\left((1-ie^{i\theta x})/(1+ie^{i\theta x})\right)+\ln^2\left((1-ie^{-i\theta x})/(1+ie^{-i\theta x})\right)}{(2^2-x^2)(4^2-x^2)}dx\\
&=\frac{-1}{4}PV\int_0^{\infty}\frac{2\ \mathrm{Re}\left(\ln^2((1+\sin(\theta x)-i\cos(\theta x))/(1-\sin(\theta x)+i\cos(\theta x)))\right)}{(2^2-x^2)(4^2-x^2)}dx\\
&=\frac{-1}{2}PV\int_0^{\infty}\frac{\mathrm{Re}\,(\ln(1+\sin(\theta x)-i\cos(\theta x))-\ln(1-\sin(\theta x)+i\cos(\theta x)))^2}{(2^2-x^2)(4^2-x^2)}dx\\
&=\frac{-1}{2}PV\int_0^{\infty}\frac{\mathrm{Re}\,(\pm i(\pi/2)-\ln|\tan((\theta x/2)-(\pi/4))|)^2}{(2^2-x^2)(4^2-x^2)}dx\\
&=\frac{-1}{2}PV\int_0^{\infty}\frac{\mathrm{Re}\left(-(\pi^2/4)+\ln^2|\tan((\theta x/2)-(\pi/4))|\pm i\pi\ln|\tan((\theta x/2)-(\pi/4))|\right)}{(2^2-x^2)(4^2-x^2)}dx\\
&=\frac{-1}{2}PV\int_0^{\infty}\frac{-(\pi^2/4)+\ln^2|\tan((\theta x/2)-(\pi/4))|}{(2^2-x^2)(4^2-x^2)}dx\\
&=\frac{-\pi}{48\,i}\left(2\left(\tan^{-1}e^{2i\theta}\right)^2-\left(\tan^{-1}e^{-2i\theta}\right)^2-\left(\left(\tan^{-1}e^{2i\theta}\right)^2-\left(\tan^{-1}e^{-2i\theta}\right)^2\right)\right).
\end{aligned} \quad (87)$$

$$\begin{aligned}
&\therefore PV\int_0^{\infty}\frac{\ln^2|\tan((\theta x/2)-(\pi/4))|}{(2^2-x^2)(4^2-x^2)}dx\\
&=-\frac{\pi}{48\,i}\left(2\left(\tan^{-1}e^{2i\theta}\right)^2-\left(\tan^{-1}e^{-2i\theta}\right)^2-\left(\left(\tan^{-1}e^{2i\theta}\right)^2-\left(\tan^{-1}e^{-2i\theta}\right)^2\right)\right).
\end{aligned}$$

## 5. Conclusion

In this study, we present new general theorems to simplify the calculation of improper integrals with principal values. These outcomes can establish many formulas of improper integrals and solve them directly without requiring complicated computations or computer softwares. We introduce some remarks to explain and analyze our results.

(i) The introduced theorems are considered powerful tools in generating and solving improper integrals and checking the accuracy of the results obtained by using other technical methods or approximating in solving similar examples

(ii) These results can be presented in tables of integrals with various values of functions and new results

(iii) The integrals obtained in this study cannot be solved manually (simply) or by computer software such as Mathematica and Maple. In future work, we will generalize the proposed results and theorems and implement them to make new tables of integrals and introduce more applications. Additionally, these results can be used to solve integral and differential equations

## Data Availability

No data were used to support this study.

## Conflicts of Interest

The authors declare that they have no conflicts of interest.